\documentclass{amsart}
\usepackage{amssymb}

\usepackage{graphicx}

\newtheorem{theorem}{Theorem}[section]
\theoremstyle{plain} \numberwithin{equation}{section}

\newtheorem{definition}[theorem]{Definition}
\newtheorem{example}[theorem]{Example}

\newtheorem{remark}[theorem]{Remark}

\begin{document}
\title{Numerical simulation of optimal transport paths}
\author{Qinglan Xia}
\address{University of California at Davis\\
Department of Mathematics\\
Davis,CA,95616} \email{qlxia@math.ucdavis.edu}
\urladdr{http://math.ucdavis.edu/\symbol{126}qlxia}
\thanks{This work is supported by an NSF grant DMS-0710714.}

\begin{abstract}
This article provides numerical simulation of an optimal transport path from
a single source to an atomic measure of equal total mass. We first construct
an initial transport path, and then modify the path as much as possible by
using both local and global minimization algorithms.
\end{abstract}

\maketitle

\section{Introduction}

This article aims at providing numerical simulations of optimal transport
paths studied earlier in \cite{xia1},\cite{xia2},\cite{xia3} \cite{xia4}
etc. The theory of optimal transport paths was motivated by studying a
phenomenon in optimal transportation where a transport system with a
branching structure may be more cost efficient than the one with a
``linear'' structure. Trees, railways, lightning, electric power supply, the
circulatory system, the river channel networks, and cardiovascular systems
are some common examples. We used the concept of optimal transport paths
between probability measures to model such transport systems in \cite{xia1}.
Some related works are \cite{buttazzo},\cite{DH},\cite{gilbert},\cite{msm}
and others. Also, in \cite{xia5}, we showed that an optimal transport path
is exactly a geodesic in the sense of metric geometry on the metric space of
probability measures with a suitable metric. As a result, people are
interested in knowing what an optimal transport path look like numerically.
This is the motivation of this article. Currently, we are using optimal
transport paths generated here to model blood vessel structures found in
placentas of human babies and also river channel networks. More application
of optimal transport paths is expected in modeling systems simulating what
we seen in the nature.

This article is organized as follows. After briefly recalling some
preliminary definitions about optimal transport paths, we study algorithms
for generating optimal transport paths. The idea is as follows. We first
consider how to generate an initial transport path, then study how to reduce
the cost by modifying the initial transport path as much as possible. We not
only use an algorithm of local minimization but also a global one. Topology
of transport paths are not preserved during the modification process. These
algorithms do not necessarily provide us a perfect optimal transport path,
but it approximate an optimal transport path very well. Using eyes of a
human being, the author cannot observe a better transport path than the path
generated here.

\section{Preliminaries}

We first recall some concepts about optimal transport paths between measures
as studied in \cite{xia1}. Let $X$ be a convex compact subset in $\mathbb{R}%
^{d}$. For any $x\in X$, let $\delta _{x}$ be the Dirac measure centered at $%
x$. An atomic measure in $X$ is in the form of
\begin{equation*}
\sum_{i=1}^{k}a_{i}\delta _{x_{i}}
\end{equation*}%
with distinct points $x_{i}\in X$, and $a_{i}>0$ for each $i=1,\cdots, k$.
Let $\mathbf{a}$ and $\mathbf{b}$ be two fixed atomic measures in the form
of
\begin{equation}
\mathbf{a}=\sum_{i=1}^{k}a_{i}\delta _{x_{i}}\text{ and }\mathbf{b}%
=\sum_{j=1}^{l}b_{j}\delta _{y_{j}}  \label{atomic_measures}
\end{equation}%
of equal total mass
\begin{equation*}
\sum_{i=1}^{k}a_{i}=\sum_{j=1}^{l}b_{j}.
\end{equation*}

\begin{definition}
A \textit{transport path} from $\mathbf{a}$ to $\mathbf{b}$ is a weighted
directed graph $G$ consists of a vertex set $V\left( G\right) $, a directed
edge set $E\left( G\right) $ and a weight function
\begin{equation*}
w:E\left( G\right) \rightarrow \left( 0,+\infty \right)
\end{equation*}%
such that $\left\{ x_{1,}x_{2,\cdots ,}x_{k}\right\} \cup \left\{
y_{1},y_{2},\cdots ,y_{l}\right\} \subset V\left( G\right) $ and for any
vertex $v\in V\left( G\right) ,$%
\begin{equation}
\sum_{\substack{ e\in E\left( G\right)  \\ e^{-}=v}}w\left( e\right) =\sum
_{\substack{ e\in E\left( G\right)  \\ e^{+}=v}}w\left( e\right) +\left\{
\begin{array}{cc}
a_{i}, & \text{if }v=x_{i}\text{ for some }i=1,\cdots ,k \\
-b_{j}, & \text{if }v=y_{j}\text{ for some }j=1,\cdots ,l \\
0, & \text{otherwise}%
\end{array}%
\right.  \label{balance_equation}
\end{equation}%
where $e^{-}$ and $e^{+}$denotes the starting and ending endpoints of each
directed edge $e\in E\left( G\right) $.
\end{definition}

\begin{remark}
The balance equation (\ref{balance_equation}) simply means that the total
mass flows into $v$ equals to the total mass flows out of $v$. When $G$ is
viewed as a polyhedral chain or current, (\ref{balance_equation}) can be
simply expressed as
\begin{equation*}
\partial G=\mathbf{b}-\mathbf{a}\text{.}
\end{equation*}
\end{remark}

Let $Path(\mathbf{a},\mathbf{b})$ be the space of all transport paths from $%
\mathbf{a}$ to $\mathbf{b}$.

\begin{definition}
For any $\alpha \leq 1$, and any $G\in Path(\mathbf{a},\mathbf{b})$, define
\begin{equation*}
\mathbf{M}_{\alpha }\left( G\right) :=\sum_{e\in E\left( G\right) }\left[
w\left( e\right) \right] ^{\alpha }length\left( e\right) .
\end{equation*}
\end{definition}

In \cite[Proposition 2.1]{xia1}, we showed that for any transport path $G\in
Path\left( \mathbf{a,b}\right) $, there exists another transport path $%
\tilde{G}\in Path\left( \mathbf{a,b}\right) $ such that
\begin{equation*}
\mathbf{M}_{\alpha }\left( \tilde{G}\right) \leq \mathbf{M}_{\alpha }\left(
G\right) ,
\end{equation*}%
vertices $V\left( \tilde{G}\right) \subset V\left( G\right) $ and $\tilde{G}$
contains no cycles. Here, a weighted directed graph $G=\left\{ V\left(
G\right) ,E\left( G\right) ,w:E\left( G\right) \rightarrow (0,1]\right\} $
contains a \textit{cycle} if for some $k\geq 3$, there exists a list of
distinct vertices $\left\{ v_{1},v_{2},\cdots ,v_{k}\right\} $ in $V\left(
G\right) $ such that for each $i=1,\cdots ,k$, either\ the segment $\left[
v_{i},v_{i+1}\right] $ or $\left[ v_{i+1},v_{i}\right] $ is a directed edge
in $E(G)$, with the agreement that $v_{k+1}=v_{1}$. When a directed graph $G$
contains no cycles, it becomes a directed tree.

An $\mathbf{M}_{\alpha}$ minimizer in $Path(\mathbf{a},\mathbf{b})$ is
called an \textit{optimal transport path} from $\mathbf{a},\mathbf{b}$.

\section{Simulation of optimal transport paths from a single source}

Let $X$ be a convex subset in $\mathbb{R}^{d}$. Given two atomic measures in
the form of
\begin{equation}
\mathbf{a}=m\delta _{O}\text{ and }\mathbf{b}=\sum_{i=1}^{N}m_{i}\delta
_{y_{i}}\text{ with }m=\sum_{i=1}^{N}m_{i}  \label{atomic_delta}
\end{equation}%
in $X$ of equal total mass, we are interested in seeing what an optimal
transport path $G$ from the single source $\mathbf{a}$ to $\mathbf{b}$ look
like numerically. 
%
%

If $N=1$, then $G$ is clearly consisting of only one edge $\left[ O,y_{1}%
\right] $ with weight $m$. If $N=2$, then we can calculate the optimal
transport path as follows.

\subsection{One source to two targets}

Suppose there are two atomic measures
\begin{equation}
\mu =m_{O}\delta _{O}\text{ and }\nu =m_{P}\delta _{P}+m_{Q}\delta _{Q},
\label{two2one}
\end{equation}%
with $m_{O}=m_{P}+m_{Q}$ for three points $P,Q,O$ in the space $X$. To find
an optimal transport path from $\mu $ to $\nu $, we need to minimize the
function%
\begin{equation*}
f\left( B\right) =\left( m_{O}\right) ^{\alpha }\left| \overrightarrow{OB}%
\right| +\left( m_{P}\right) ^{\alpha }\left| \overrightarrow{BP}\right|
+\left( m_{Q}\right) ^{\alpha }\left| \overrightarrow{BQ}\right|
\end{equation*}%
among all points $B$ in the triangle $\triangle POQ$. Here, we use the
notation $\overrightarrow{BP}$ etc. to denote the vector $P-B$ in $X$, and
let $\left| \overrightarrow{BP}\right| $ be the magnitude of this vector.
Since $f\left( B\right) $ is a continuous function on a compact set, $\ f$
must achieve its minimum at some point $B^{\ast }$. Indeed, we can find $%
B^{\ast }$ as follows. Suppose $B^{\ast }$ is located in the interior of the
triangle $\triangle POQ$, then it must satisfy the balance equation
\begin{equation*}
\left( m_{O}\right) ^{\alpha }\frac{\overrightarrow{OB}}{\left|
\overrightarrow{OB}\right| }+\left( m_{P}\right) ^{\alpha }\frac{%
\overrightarrow{BP}}{\left| \overrightarrow{BP}\right| }+\left( m_{Q}\right)
^{\alpha }\frac{\overrightarrow{BQ}}{\left| \overrightarrow{BQ}\right| }=%
\vec{0}
\end{equation*}%
at $B=B^{\ast }$. From it, one can easily find the angles
\begin{equation*}
\measuredangle OB^{\ast }P=\theta _{1},\measuredangle OB^{\ast }Q=\theta _{2}%
\text{ and }\measuredangle PB^{\ast }Q=\theta _{3}
\end{equation*}%
where
\begin{equation}
\theta _{1}=\cos ^{-1}\left( \frac{k_{2}-k_{1}-1}{2\sqrt{k_{1}}}\right)
,\theta _{2}=\cos ^{-1}\left( \frac{k_{1}-k_{2}-1}{2\sqrt{k_{2}}}\right)
,\theta _{3}=\cos ^{-1}\left( \frac{1-k_{1}-k_{2}}{2\sqrt{k_{1}k_{2}}}\right)
\label{angles}
\end{equation}%
for
\begin{equation*}
k_{1}=\left( \frac{m_{P}}{m_{O}}\right) ^{2\alpha },k_{2}=\left( \frac{m_{Q}%
}{m_{O}}\right) ^{2\alpha }\text{.}
\end{equation*}

Let $M$ (and $H$) be the projection of the point $Q$ (and $P$, respectively)
along the segment $\overrightarrow{OP}$ (and $\overrightarrow{OQ}$
respectively). Then, the centers $R$ (, and $S$) of the circles passing
through the triangles $\triangle OB^{\ast }P$ (and $\triangle OB^{\ast }Q$
respectively) is given by

\begin{eqnarray*}
R &=&\frac{O+P}{2}-\frac{\cot \theta _{1}}{2}\frac{\overrightarrow{QM}}{%
\left| \overrightarrow{QM}\right| }\left| \overrightarrow{OP}\right| \\
S &=&\frac{O+Q}{2}-\frac{\cot \theta _{2}}{2}\frac{\overrightarrow{PH}}{%
\left| \overrightarrow{PH}\right| }\left| \overrightarrow{OQ}\right|
\end{eqnarray*}%
where
\begin{equation*}
\overrightarrow{QM}=\frac{\overrightarrow{OP}\cdot \overrightarrow{OQ}}{%
\left| \overrightarrow{OP}\right| ^{2}}\overrightarrow{OP}-\overrightarrow{OQ%
}\text{ and }\overrightarrow{PH}=\frac{\overrightarrow{OP}\cdot
\overrightarrow{OQ}}{\left| \overrightarrow{OQ}\right| ^{2}}\overrightarrow{%
OQ}-\overrightarrow{OP}.
\end{equation*}
By (\ref{angles}),
\begin{equation*}
\cot \theta _{1}=\frac{k_{2}-k_{1}-1}{\sqrt{4k_{1}-\left(
k_{2}-k_{1}-1\right) ^{2}}}\text{ and }\cot \theta _{2}=\frac{k_{1}-k_{2}-1}{%
\sqrt{4k_{2}-\left( k_{1}-k_{2}-1\right) ^{2}}}.
\end{equation*}

Now, $B^{\ast }$ is just the reflection of the point $O$ along the segment $%
RS$. That is,

\begin{equation*}
B^{\ast }=2\left[ (1-\lambda )R+\lambda S\right] -O\text{ \ with }\lambda =%
\frac{\overrightarrow{RO}\cdot \overrightarrow{RS}}{\left| \overrightarrow{RS%
}\right| ^{2}}
\end{equation*}%
whenever $B^{\ast }$ is located in the interior of the triangle $\triangle
POQ$. Note that a necessary condition for $B^{\ast }$ being located in the
interior of the triangle $\triangle POQ$ is the angles must satisfy
\begin{equation*}
\measuredangle OQP<\theta _{1}\text{,}\measuredangle OPQ<\theta _{2}\text{
and }\measuredangle POQ<\theta _{3}\text{. }
\end{equation*}%
In case the condition fails, we have three degenerate cases. If the angle $%
\measuredangle POQ\geq \theta _{3}$, then take $B^{\ast }$ to be $O$ and we
get a ``V-shaped'' path. If the angle $\measuredangle OQP\geq \theta _{1}$
and $\measuredangle POQ<\theta _{3}$, then take $B^{\ast }$ to be $Q$. If
the angle $\measuredangle OPQ\geq \theta _{2}$ and $\measuredangle
POQ<\theta _{3}$, then\ take $B^{\ast }$ to be $P$.
\begin{figure}[h]
\includegraphics[height=2in]{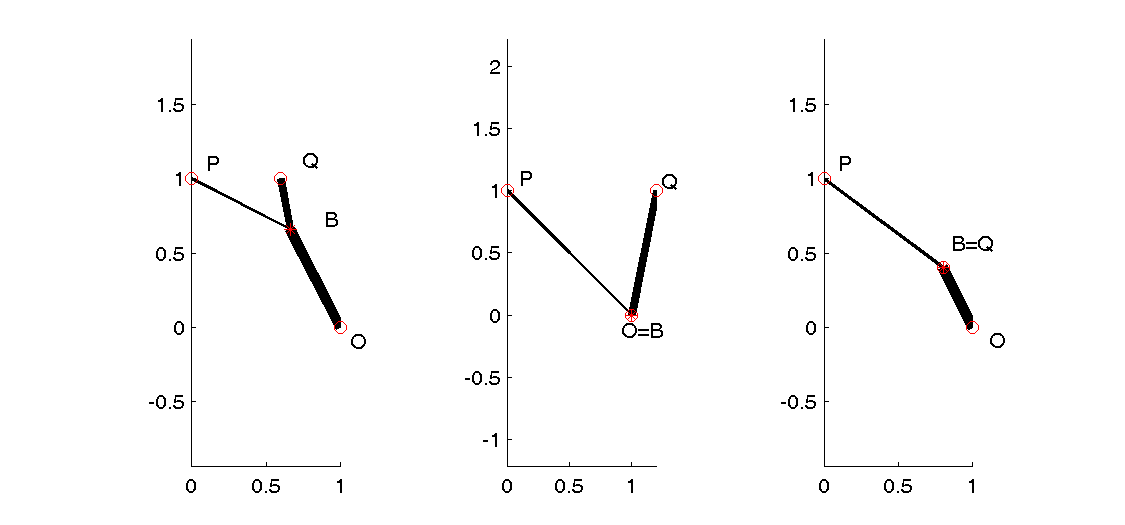}
\caption{Cases for transporting $\protect\mu$ to $\protect\nu$}
\end{figure}

As a result, given $\mu $ and $\nu $ in (\ref{two2one}), we achieved a
formula for finding $B^{\ast }$. The optimal transport path $G$ from $\mu $
to $\nu $ has at most three edges: $\left[ B^{\ast },P\right] $ with weight $%
m_{P}$, $\left[ B^{\ast },Q\right] $ with weight $m_{Q}$ and $\left[
O,B^{\ast }\right] $ with weight $m_{O}$.

We denote the point $B^{\ast }$ by $V\left( \mu ,\nu \right) $ and let
\begin{equation}
g\left( \mu ,\nu \right) =f\left( O\right) -f\left( B^{\ast }\right)
\label{advantage}
\end{equation}%
which gives the advantage of taking a ``Y-shaped'' path over taking a
``V-shaped'' path.

\subsection{The construction of an initial transport path}

When $N\geq 3$, it is not that easy (or may be impossible) to find an exact
solution of an optimal transport path. Instead, we would like to find an
approximately optimal transport path. The idea is to construct an initial
transport path $G\in Path\left( \mathbf{a,b}\right) $ and then modify $G$ as
much as possible until we can not reduce the cost of $G$ any further.
%

\subsubsection{The method of transporting small number of points\label%
{small_number}}

If $N\leq 2$, then we have found the optimal transport path $G$ as above. If
$2<N\leq K$ for a small given number $K$, then for any pair $1\leq i<j\leq N$%
, let
\begin{equation*}
g_{ij}=g\left( \left( m_{i}+m_{j}\right) \delta _{O},m_{i}\delta
_{y_{i}}+m_{j}\delta _{y_{j}}\right)
\end{equation*}%
where the function $g$ is defined as in (\ref{advantage}). Suppose the
maximum of $\left\{ g_{ij}\right\} $ is achieved at $1\leq i^{\ast }<j^{\ast
}\leq N$. Then, the desired path $G$ is given recursively by
\begin{equation*}
G=\tilde{G}+m_{i^{\ast }}\left[ B^{\ast },y_{i^{\ast }}\right] +m_{j^{\ast
}}[B^{\ast },y_{j^{\ast }}],
\end{equation*}%
where $B^{\ast }=V\left( \left( m_{i^{\ast }}+m_{j^{\ast }}\right) \delta
_{O},m_{i^{\ast }}\delta _{y_{i^{\ast }}}+m_{j^{\ast }}\delta _{y_{j^{\ast
}}}\right) $ is the point in $X$ given by (\ref{advantage}), and $\tilde{G}$
is the path from $\mathbf{a}$ to $\mathbf{\tilde{b}}=\mathbf{b}-m_{i^{\ast
}}\delta _{y_{i^{\ast }}}-m_{j^{\ast }}\delta _{y_{j^{\ast }}}+\left(
m_{i^{\ast }}+m_{j^{\ast }}\right) \delta _{B^{\ast }} $ achieved by
recursively applying this algorithm.

\subsubsection{The subdivision method}

To construct an initial transport path in $Path\left( \mathbf{a,b}\right) $,
one may simply take a trivial transport path
\begin{equation*}
\sum_{i=1}^{N}m_{i}\left[ O,y_{i}\right] .
\end{equation*}%
This is an allowable transport path in $Path\left( \mathbf{a,b}\right) $.
Nevertheless, the degree of the vertex $O$ (i.e. the total number of edges
in $G$ having $O$ as an endpoint) is $N$. Then, it might become time
consuming later for modifying the path at the vertex $O$ when $N$ is very
large. Instead, we use the following subdivision method to construct an
initial transport path $G_{sd}\left( \mathbf{a},\mathbf{b}\right) $, which
contains no cycles and has degree at most $K$ at every vertex for some given
$K$ defined below.

Let
\begin{equation*}
\lambda =\left\{
\begin{array}{cc}
3, & \text{if }d=2 \\
2, & \text{if }d\geq 3%
\end{array}%
\right.
\end{equation*}%
and $K=\lambda ^{d}$ where $d$ is the dimension of the ambient space $%
\mathbb{R}^{d}$.

\textbf{Algorithm} (subdivision method):

\textit{Input: two atomic measures $\mathbf{a},\mathbf{b}$ in the form of (%
\ref{atomic_delta}) and a parameter $\alpha \leq 1$};

\textit{Output: a transport path $G\in Path\left( \mathbf{a,b}\right) $ with
degree$\left( v\right) \leq K$ for each $v\in V(G)$.}

If $N\leq K$, then we use the method of transporting small number of points\
described above to construct a transport path from $\mathbf{a}$ to $\mathbf{b%
}$.

If $N>K$, then let $Q$ be a cube in $R^{d}$ that contains the supports of
both $a$ and $b\,$. We may split the cube $Q$ into totally $K=\lambda^d$
smaller cubes $\left\{ Q_{i}\right\} _{i=1}^{K}$ of size equal to $\frac{1}{%
\lambda }$ of the size of $Q$. \ For each $i=1,\cdots ,K$, let $G_{i}$ be
the path $G_{sd}\left( \mathbf{b}\left( Q_{i}\right) \delta _{c\left(
Q_{i}\right) },\mathbf{b}\lfloor _{Q_{i}}\right) $ from the center $c\left(
Q_{i}\right) $ of the smaller cube $Q_{i}$ to the restriction of $b$ in $%
Q_{i}$ achieved by recursively applying this algorithm. Also, let $G_{0}$ be
the path from $\mathbf{a}$ to $\sum_{i=1}^{K}\mathbf{b}\left( Q_{i}\right)
\delta _{c_{\left( Q_{i}\right) }}$ by using the method of transporting
small number of points. Then,
\begin{equation*}
G=\sum_{i=0}^{K}G_{i}
\end{equation*}%
provides the desired path $G_{sd}\left( \mathbf{a},\mathbf{b}\right) $ from $%
\mathbf{a}$ to $\mathbf{b}$.

\begin{figure}[h]
\includegraphics[height=5in]{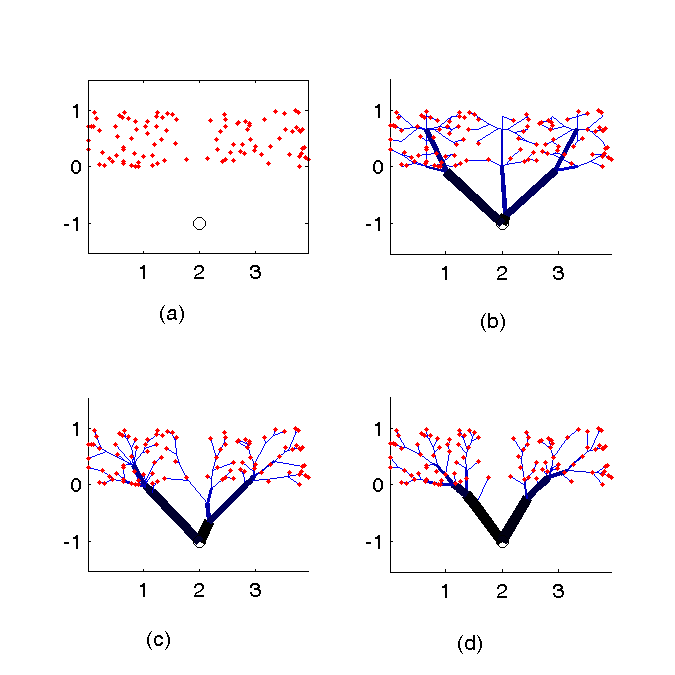} %
\caption{(a) A single source $(2,-1)$ and the targeting measure represented
by 100 random points. (b) An initial transport path constructed by the
subdivision method. (c)A modified transport path achieved by repeated
modifying the initial path in (b) using the local minimization. (d)An
optimal transport path achieved by modifying the path in (c) using the
global minimization method.}
\end{figure}

\subsection{Modification of an existing transport path}

Now, suppose $G$ is an existing transport path from $\mathbf{a}$ to $\mathbf{%
b}$ that contains no cycles. We want to modify $G$ to reduce the transport
cost as much as possible. Before describing algorithms, we first introduce
some concepts about vertices of an transport path $G$.

For any two vertices $v,u\in V\left( G\right) $, we say that $v$ is an
\textit{ancestor} of $u$ and $u$ is a \textit{descendant} of $v$, if there
exists a list of vertices $v_{1}=v,v_{2},\cdots ,v_{h-1},v_{h}=u$ such that
each $\left[ v_{i},v_{i+1}\right] $ is a directed edge in $E\left( G\right) $
for $i=1,\cdots ,h-1$. Also, if $\left[ v,u\right] $ is a directed edge in $%
E\left( G\right) $, then we say that $v$ is a \textit{parent} of $u$ and $u$
is a \textit{child} of $v$.

For each vertex $u\in V\left( G\right) \backslash \left\{ O\right\} $, $u$
has exactly one parent $p\left( u\right) \in V\left( G\right) $ because $G$
contains no cycles and has a single source $\mathbf{a}=m\delta _{O}$. Let $%
m\left( u\right) $ be the associated weight on the directed edge $\left[
p\left( u\right) ,u\right] $ in $E\left( G\right) $ for each $u\in
V(G)\setminus \left\{ O\right\}$, and also set $m\left( O\right) =m$. Note
that $m\left( v\right) \geq m\left( u\right) $ whenever $v $ is an ancestor
of $u$. Moreover, the vertex $O$ is always an ancestor of each $u$. That is,
there exists a list of vertices $v_{1}$,$v_{2}$,$\cdots ,v_{k}$ in $V\left(
G\right) $ such that $\left[ v_{i},v_{i+1}\right] \in E\left( G\right) $
with $v_{1}=O$ and $v_{k}=u$. Then, for each $t\in \left[ -m\left( u\right)
,m\left( u\right) \right] $, we consider the path
\begin{equation*}
R\left( G;t,u\right) :=G-\sum_{i=1}^{k-1}t\left[ v_{i},v_{i+1}\right] \in
Path\left( \mathbf{a}-t\delta _{O}+t\delta _{u},\mathbf{b}\right) .
\end{equation*}
When $t>0$, we say that a mass of $t$ is \textit{removed} from the path $G$
at vertex $u$, When $t<0$, we say that a mass of $t$ is \textit{added} to
the path $G$ at vertex $u$. Moreover, the potential function of $G$ at a
vertex $u\in V\left( G\right) $ is defined by
\begin{equation}
P_{G}\left( u,t\right) =\left\{
\begin{array}{cc}
P_{G}\left( p\left( u\right) ,t\right) +\left| p\left( u\right) -u\right|
\left[ m\left( u\right) ^{\alpha }-\left( m\left( u\right) -t\right)
^{\alpha }\right] , & u\neq O \\
0, & u=O%
\end{array}%
\right.  \label{potential}
\end{equation}

for $t\in \left[ -m\left( u\right) ,m\left( u\right) \right] $. 
Note that $P_{G}\left( u,t\right) $ has the same sign as $t$.

\subsubsection{local minimization}

We first use a local minimization method to modify any existing transport
path $G$ containing no cycles.

\textit{Input: a transport path $G\in Path\left( \mathbf{a,b}\right) $
containing no cycles and $\alpha <1 $};

\textit{Output: a locally optimized path $\tilde{G}\in Path\left( \mathbf{a,b%
}\right) $ with $M_{\alpha }\left( \tilde{G}\right) \leq M_{\alpha }\left(
G\right) $}.

\textit{Idea: For each vertex }$u$\textit{\ in }$G$\textit{, replace }$%
G_{old}\left( u\right) $\textit{\ by }$G_{new}\left( u\right) $\textit{\
whenever }$M_{\alpha }\left( G_{old}\left( u\right) \right) >M_{\alpha
}\left( G_{new}\left( u\right) \right) $\textit{.}

\bigskip Here, for each vertex $u$ of $G$, two transport paths $%
G_{old}\left( u\right) $ and $G_{new}\left( u\right) $ are defined as
follows. Let
\begin{equation*}
\mu _{C}=\sum_{h\in V\left( G\right) ,p\left( h\right) =u}m\left( h\right)
\delta _{h}\text{ and }\mu _{P}=m\left( u\right) \delta _{p\left( u\right) }
\end{equation*}%
be two atomic measures corresponding to the children and the parent of $u$.
Then,
\begin{equation*}
G_{old}\left( u\right) =\sum_{h\in V\left( G\right) ,p\left( h\right)
=u}m\left( h\right) \left[ u,h\right] +m\left( u\right) \left[ p\left(
u\right) ,u\right] \in Path\left( \mu _{P},\mu _{C}\right)
\end{equation*}%
is the union of all weighted edges in $G$ sharing $u$ as their common
endpoint. On the other hand, one may generate another path $G_{new}\left(
u\right) \in Path\left( \mu _{P},\mu _{C}\right) $ by using the method of
transporting small number of points stated in \ref{small_number}.

If
\begin{equation*}
M_{\alpha }\left( G_{old}\left( u\right) \right) >M_{\alpha }\left(
G_{new}\left( u\right) \right) \text{,}
\end{equation*}%
then by replacing $G_{old}\left( u\right) $ by $G_{new}\left( u\right) $ in $%
G$, we get a new path
\begin{equation*}
\tilde{G}=G-G_{old}\left( u\right) +G_{new}\left( u\right) \in Path\left(
\mathbf{a,b}\right)
\end{equation*}%
and $M_{\alpha }\left( \tilde{G}\right) \leq M_{\alpha }\left( G\right)
-M_{\alpha }\left( G_{old}\left( u\right) \right) +M_{\alpha }\left(
G_{new}\left( u\right) \right) <M_{\alpha }\left( G\right) $. So, $\tilde{G%
}$ is a transport path with less cost. Replace $G$ by this modified
path $\tilde{G}$, and continue this process for all vertices of $G$
until one can not reduce the cost any further.

The main drawback of this algorithm is that the result is only local
minimization rather than global minimization. For instance, edges may
intersect with each other. Sometimes, using eyes of a human being, one can
easily observe a better transport path. To overcome these drawbacks, we
adopt the following algorithm.

\subsubsection{global minimization}

Now, we introduce the following algorithm of global minimization:

\textit{Input: two probability measures }$\mathbf{a}$\textit{, }$\mathbf{b}$%
\textit{\ in the form of (\ref{atomic_delta}) and a parameter }$\alpha \leq
1 $;

\textit{Output: an approximately $M_{\alpha }$ optimal transport path $G\in
Path\left( \mathbf{a,b}\right) $.}

step 1: construct a transport path $G$ from $\mathbf{a}$ to $
\mathbf{b}$ using the subdivision method;

step 2: modify the existing path $G$ using the local minimization
method;

step 3: subdivide long edges of $G$ into shorter edges;

step 4: for each vertex $u$ of $G$, remove a mass of $m\left( u\right) $ at vertex $%
u$ from the path $G$; change the parent $p(u)$ of $u$ to a better
one if possible and then add back a mass of $m\left( u\right) $ at
vertex $u$. More precisely,

substep 1: A list of potential parents of $u$ is defined as
\begin{equation*}
L\left( u\right) =\left\{ v\in V\left( G\right) :\left| v%
-u\right| \leq \sigma \text{, and }v\text{ is not a descendant of }%
u\right\} \text{,}
\end{equation*}

\bigskip where
\begin{equation*}
\sigma =\frac{P_{G}\left( u,m\left( u\right) \right) }{\left[ m\left(
u\right) \right] ^{\alpha }}
\end{equation*}%
and $P_{G}$ is defined in (\ref{potential}). Note that the parent
$p(u)$ is automatically in $L\left( u\right) $ because
\begin{equation*}
\sigma =\frac{P_{G}\left( p(u),m\left( u\right) \right) +\left|
p(u)-u\right| m\left( u\right) ^{\alpha }}{\left[ m\left( u\right)
\right] ^{\alpha }}\geq \left| p(u)-u\right|.
\end{equation*}%

substep 2: By removing a mass of $m\left( u\right) $ at vertex $u$ from the
path $G$, we get another path
\begin{equation*}
\tilde{G}=R\left( G;m\left( u\right) ,u\right) .
\end{equation*}

substep 3: For each $v\in L\left( u\right) \setminus \left\{
p(u)\right\} $, let%
\begin{equation*}
c\left( v\right) =-P_{\tilde{G}}\left( v,-m\left( u\right) \right) ,
\end{equation*}%
where $P_{\tilde{G}}$ is defined as in (\ref{potential}) with $G$
replaced by $\tilde{G}$. The number $c\left( v\right) $ measures the
extra cost of transporting a mass of $m\left( u\right) $ on the
system $\tilde{G}$ from the source $O$ to the vertex $u$ via the
vertex $v$.

substep 4: Find the maximum of $c\left( v\right) $ over all $v \in
L\left( u\right) \setminus \left\{ p(u)\right\} $. If $\max c\left(
v\right) >\sigma \left[ m\left( u\right) \right] ^{\alpha }$, then
we find
a better parent for the vertex $u$. In this case, suppose the maximum of $%
c\left( v\right) $ is achieved at $v^{\ast }$. Then, let
\begin{equation*}
G^{\ast }=R\left( \tilde{G};-m\left( u\right) ,v^{\ast }\right) +m\left(
u\right) \left[ v^{\ast },u\right] .
\end{equation*}%
That is, we change the parent of $u$ from $p(u)$ to $v^{\ast }$ and then add a mass $%
m\left( u\right) $ at the vertex $u$ to the modified path. For convenience,
we still denote the final modified transport path $G^{\ast }$ by $G$.

step 5: Repeat steps 2-4 until one can not reduce the cost any further.

\section{Examples}

\begin{example}
Let $\left\{ y_{i}\right\} $\bigskip\ be 50 random points in the square $%
\left[ 0,1\right] \times \left[ 0,1\right] $. Then, $\left\{ y_{i}\right\} $
determines an atomic probability measure
\begin{equation*}
\mathbf{b}=\sum_{i=1}^{50}\frac{1}{50}\delta _{y_{i}}.
\end{equation*}%
Let \textbf{$a$}$=\delta _{O}$ where $O=\left( 0,0\right) $ is the origin.
Then an optimal transport path from $\mathbf{a}$ to $\mathbf{b}$ looks like
the following figures with $\alpha =1,0.75,0.5$ and $0.25$ respectively:

\includegraphics[height=1.75in]{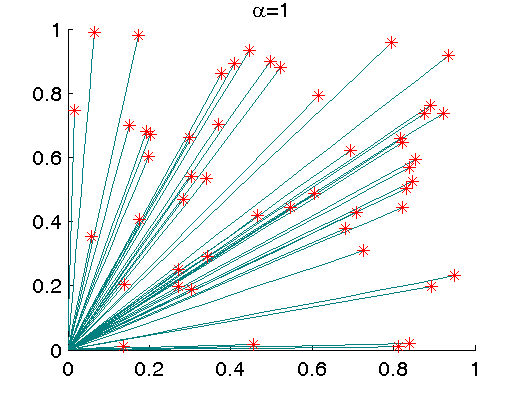}%
\includegraphics[height=1.75in]{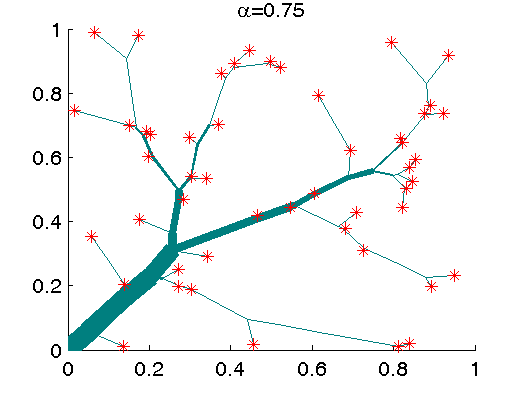}

\includegraphics[height=1.75in]{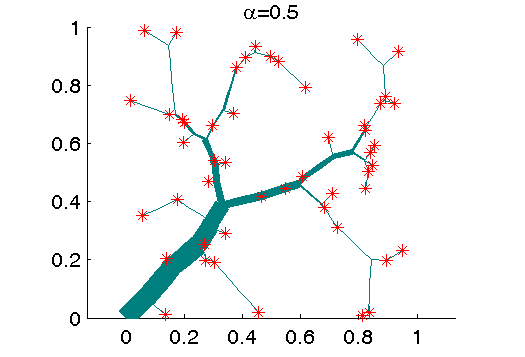}%
\includegraphics[height=1.75in]{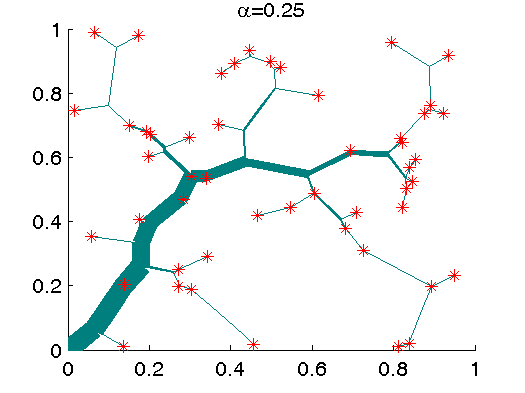}
\end{example}

\begin{example}
Let $\left\{ y_{i}\right\} $ be 100 random points in the rectangle $\left[
-2.5,2.5\right] \times \left[ 0,1\right] $. Then, $\left\{ y_{i}\right\} $
determines an atomic probability measure $\mathbf{b}=\sum_{i=1}^{100}\frac{1%
}{100}\delta _{y_{i}}.$ Let $\mathbf{a}=\delta _{O}$ where $O=\left(
0,0\right) $ is the origin, and let $\alpha =0.85$. Then an optimal
transport path from $\mathbf{a}$ to $\mathbf{b}$ looks like the following
figure.

\includegraphics[height=3.5in]{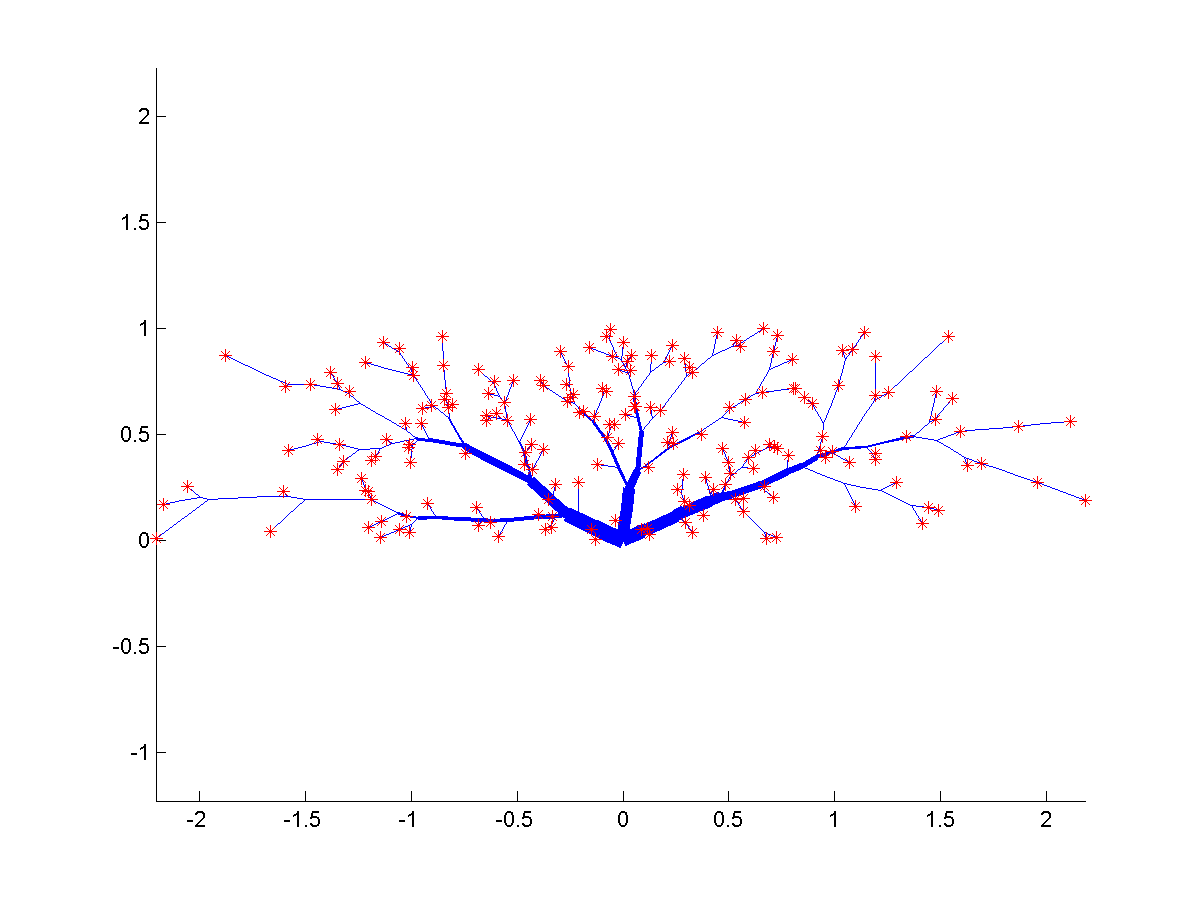}
\end{example}

\begin{example}
Optimal transport paths from the center to the unit circle. Here, the unit
circle is represented by 400 points uniformly distributed on the circle. The
parameter $\alpha =0.75$ in the first figure and $\alpha =0.95$ in the
second one.

\includegraphics[height=1.75in]{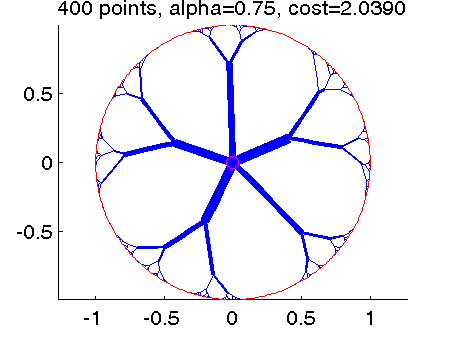} %
\includegraphics[height=1.75in]{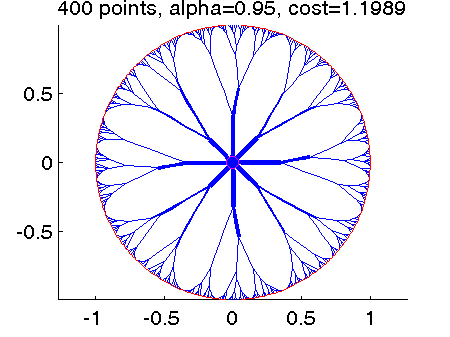}
\end{example}

\begin{example}
Optimal transport paths from the center to the unit disk. The first one is
using random generated points in the disk with $\alpha =2/3$ while the
second one use uniformly generated points in the disk with $\alpha =0.75$.

\includegraphics[height=2.5in]{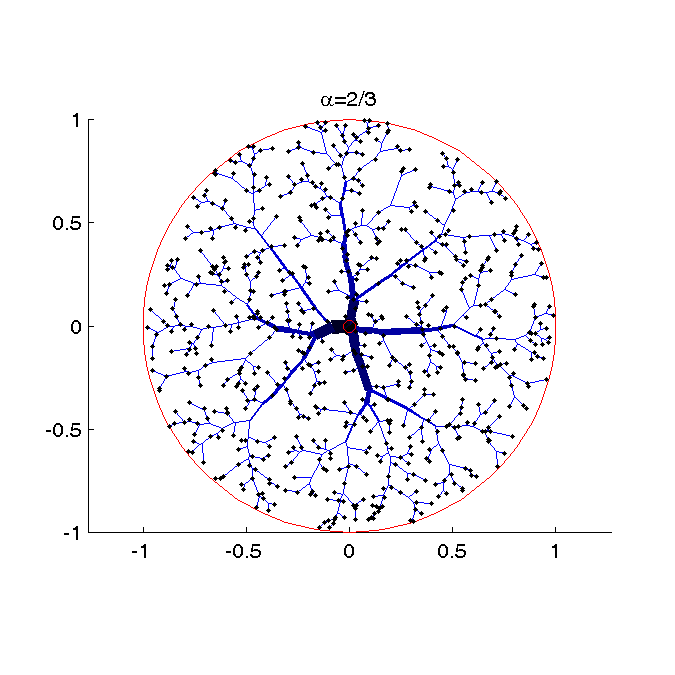} %
\includegraphics[height=2.5in]{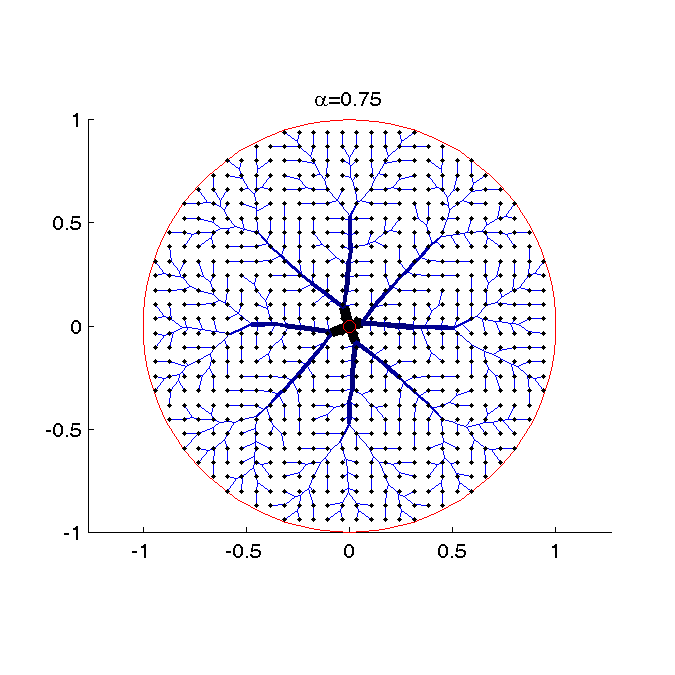}
\end{example}

\begin{example}
An optimal transport path from a point on the boundary to the unit square,
which is represented by 400 randomly generated points, with $\alpha=0.85$.

\includegraphics[height=2.2in]{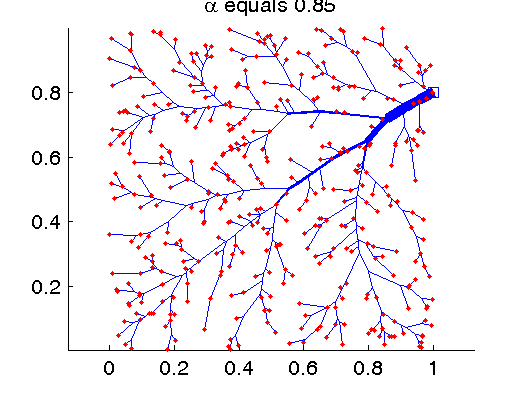}
\end{example}

\begin{example}
An optimal transport path modeling blood vessels in a placenta of new baby

\includegraphics[height=3.5in]{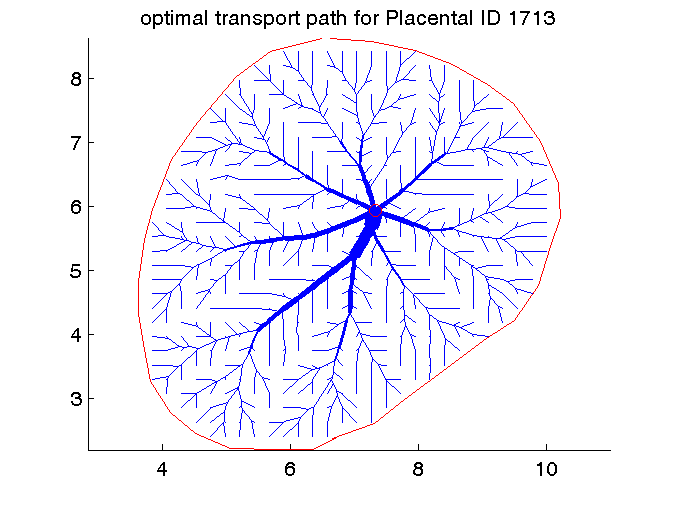}
\end{example}

\end{document}